\documentclass[10pt,letterpaper]{article}
\usepackage{blindtext,titlefoot}
\usepackage{authblk}
\usepackage{marvosym}
\usepackage[utf8]{inputenc}
\usepackage[english]{babel}
\usepackage{amsmath}
\usepackage{amsfonts}
\usepackage{amssymb}
\usepackage{makeidx}
\usepackage{graphicx}
\usepackage{lmodern}
\usepackage{kpfonts}
\usepackage{dsfont}
\usepackage{cancel}
\usepackage{amsthm} 
\usepackage[left=2cm,right=2cm,top=2cm,bottom=2cm]{geometry}
\usepackage{subcaption}
\usepackage{multirow}
\usepackage{float} 
\usepackage{url}
\usepackage{breakurl}
\usepackage[breaklinks]{hyperref}
\usepackage{multicol}
\usepackage{mathtools, cuted}
\usepackage{lipsum}
\usepackage{booktabs}
\usepackage{comment}
\usepackage{cite}

\usepackage{fancyhdr}
\fancypagestyle{firstpage}{\fancyhf{}
\setcounter{page}{1}
\rfoot{\thepage}
}

\providecommand{\nset}[1]{
\mathbb{#1}
}
\providecommand{\set}[1]{
\left\{#1\right\}
}

\providecommand{\norm}[1]{
\left\lVert #1 \right\rVert
}

\providecommand{\ds}[1]{
\displaystyle #1
}
\providecommand{\der}[3]{
\dfrac{#1^{#3} }{ #1 #2^{#3}}
}

\newtheorem{theorem}{ Theorem}[section]
\newtheorem{definition}[theorem]{Definition}

\newtheorem{corollary}[theorem]{Corollary}

\usepackage{enumitem} 
\setlist[itemize]{noitemsep} 

\usepackage{titlesec} 
\titleformat{\section}[block]{\large\bfseries\scshape\centering}{\thesection.}{1em}{} 
\titleformat{\subsection}[block]{\large\bfseries\scshape\centering}{\thesubsection.}{1em}{}
\titleformat{\subsubsection}[block]{\large\bfseries\scshape\centering}{\thesubsubsection.}{1em}{} 

\title{\bfseries A nonlinear system related to investment under uncertainty solved using the fractional pseudo-Newton method}

\author[,a]{A. Torres-Hernandez  \footnote{\scriptsize  E-mail : anthony.torres@ciencias.unam.mx; ORCID: 0000-0001-6496-9505}}
\affil[a]{\normalsize Department of Physics, Faculty of Science - UNAM, Mexico}

\author[,b]{F. Brambila-Paz \footnote{\scriptsize E-mail : fernandobrambila@gmail.com; ORCID: 0000-0001-7896-6460}}
\affil[b]{\normalsize Department of Mathematics, Faculty of Science - UNAM, Mexico}

\author[,c]{J. J.  Brambila \footnote{\scriptsize E-mail : jbrambilaa@yahoo.com.mx; ORCID: 0000-0001-7640-8203}}
\affil[c]{\normalsize Department of Economics, Colegio de Postgraduados, Mexico .}

\date{}

\begin{document}

\maketitle
\thispagestyle{firstpage}


\begin{abstract}

A nonlinear algebraic equation system of two variables is numerically solved, which is derived from a nonlinear algebraic equation system of four variables, that corresponds to a mathematical model related to investment under conditions of uncertainty. The theory of investment under uncertainty scenarios proposes a model to determine when a producer must expand or close, depending on his income. The system mentioned above is solved using a fractional iterative method, valid for one and several variables, that uses the properties of fractional calculus, in particular the fact that the fractional derivatives of constants are not always zero, to find solutions of nonlinear systems.

\textbf{Mathematical subject classification:} Fractional Calculus, Numerical Analysis.

\end{abstract}

\section{Introduction}

A classic problem in mathematics, which is of common interest in physics and engineering, is finding the set of zeros of a function $f:\Omega \subset \nset{R}^n \to \nset{R}^n$, that is,

\begin{eqnarray}\label{eq:1-001}
\footnotesize
\begin{array}{c}
\set{\xi \in \Omega \ : \ \norm{f(\xi)}=0},
\end{array}
\end{eqnarray}

where $\norm{ \ \cdot \ }: \nset{R}^n \to \nset{R}$ denotes any vector norm. Although finding the zeros of a function may seem like a simple problem, in general it involves solving an \textbf{algebraic equation system} as follows

\begin{eqnarray}\label{eq:1-002}
\footnotesize
\left\{
\begin{array}{c}
\left[f\right]_1(x)=0\\
\left[f\right]_2(x)=0\\
\vdots \\
\left[f\right]_n(x)=0
\end{array}\right.,
\end{eqnarray}

where $[f]_k: \nset{R}^n \to \nset{R}$ denotes the $k$-th component of the function $f$.  It should be noted that the system of equations \eqref{eq:1-002} may represent a \textbf{linear system} or a \textbf{nonlinear system}, and in general, it is necessary to use numerical methods of the iterative type to solve it.

It is necessary to mention that the iterative methods have an intrinsic problem, since if a system has $N$ solutions it is necessary to invest time in finding $N$ initial conditions, but this problem is partially solved by combining iterative methods with fractional calculus, whose result is known as \textbf{fractional iterative methods}, since these new methods have the ability to find $N$ solutions of a system using a single initial condition. In this document, a fractional iterative method that does not explicitly depend on the fractional partial derivatives of the function for which zeros are searched is presented, this characteristic makes it an ideal iterative method to solve nonlinear systems in several variables.

\section{Fixed Point Method}

Let  $\Phi:\nset{R}^n \to \nset{R}^n$ be a function. It is possible to build a sequence $\set{x_i}_{i=0}^\infty$  by defining the following iterative method

\begin{eqnarray}\label{eq:2-001}
\footnotesize
\begin{array}{c}
x_{i+1}:=\Phi(x_i),
\end{array}
\end{eqnarray}

if is true that $x_i\to \xi\in \nset{R}^n$ and if the function $\Phi$ is continuous around $\xi$, we obtain that

\begin{eqnarray}\label{eq:2-002}
\footnotesize
\begin{array}{c}
\ds \xi=\lim_{i\to \infty}x_{i+1}=\lim_{i\to \infty}\Phi(x_i)=\Phi\left(\lim_{i\to \infty}x_i \right)=\Phi(\xi),
\end{array}
\end{eqnarray}

the above result is the reason by which the method \eqref{eq:2-001} is known as the \textbf{fixed point method}. Furthermore, the function $\Phi$ is called an \textbf{iteration function}.

\subsection{Order of Convergence}

Consider the following definition \cite{plato2003concise}

\begin{definition}
Let $\Phi: \Omega \subset \nset{R}^ n \to \nset{R}^ n $ be an iteration function with a fixed point $ \xi \in \Omega $. Then the method \eqref{eq:2-001} is called  \textbf{(locally) convergent of (at least) order $ \boldsymbol{p} $} ($ p \geq 1 $), if there are exists $ \delta> 0 $  and $ C $, a non-negative constant  with $ C <1 $ if $ p = 1 $, such that for any initial value $ x_0 \in B (\xi; \delta) $ it holds that

\begin{eqnarray}\label{eq:c2.08}
\footnotesize
\begin{array}{cc}
\norm{x_{k+1}-\xi}\leq C \norm{x_k-\xi}^p, & k=0,1,2,\cdots,
\end{array}
\end{eqnarray}

where $ C $ is called convergence factor.

\end{definition}

The following corollary \cite{stoer2013, torres2020reduction} allows characterizing the order of convergence of an iteration function $ \Phi $ with its \textbf{Jacobian matrix} .

\begin{corollary}\label{cor:2-001}
Let $\Phi:\nset{R}^n \to \nset{R}^n$ be an iteration function. If $\Phi$ defines a sequence $\set{x_i}_{i=0}^\infty$ such that $x_i\to \xi$, and if the following condition is fulfilled

\begin{eqnarray}\label{eq:c2.16}
\footnotesize
\begin{array}{c}
\ds \lim_{x\to \xi}\norm{\Phi^{(1)}(x)}\neq 0,
\end{array}
\end{eqnarray}

where $\Phi^{(1)}$ denotes the Jacobian matrix of the function $\Phi$, then $\Phi$ has an order of convergence (at least) linear  in $B(\xi;\delta)$.
\end{corollary}

\section{Fractional Pseudo-Newton Method}

The interest in fractional calculus has mainly focused on the study and development of techniques to solve differential equation systems of order non-integer \cite{hilfer00,kilbas2006theory, martinez2017applications1,martinez2017applications2,torreshern2019proposal}. Over the years, iterative methods have also been developed that use the properties of fractional derivatives to solve algebraic equation systems \cite{gao2009local,fernando2017fractional,brambila2018fractional,akgul2019fractional,torreshern2020,torres2020fractional,gdawiec2020newton,torres2020approximation}. 
These methods may be called \textbf{fractional iterative methods}, which under certain conditions, may accelerate their speed of convergence with the implementation of the Aitken's method     \cite{brambila2018fractional}.

A recently proposed fractional iterative method, which is valid for one and several variables, and which has already been used in an engineering  application \cite{torres2020fractional}, is the \textbf{fractional pseudo-Newton method} \cite{torres2020approximation,torres2020reduction}, which is given by the following expression

\begin{eqnarray}\label{eq:c2.40}
\footnotesize
\begin{array}{cc}
x_{i+1}:=\Phi(\alpha, x_i)= x_i- P_{\epsilon,\beta}(x_i) f(x_i), & i=0,1,2\cdots,
\end{array}
\end{eqnarray}

with $\alpha\in\nset{R}\setminus\nset{Z}$, in particular $\alpha\in[-2,2]\setminus\nset{Z}$ \cite{torreshern2020}, where $ P_{\epsilon, \beta} (x_i) $ is a matrix evaluated in the value $ x_i $, which is given by the following expression

\begin{eqnarray}
\footnotesize
\begin{array}{c}
P_{\epsilon,\beta}(x_i):=\left([P_{\epsilon,\beta}]_{jk}(x_i)\right)=\left( \partial_k^{\beta(\alpha,[x_i]_k)}\delta_{jk}+ \epsilon\delta_{jk}  \right)_{x_i},
\end{array}
\end{eqnarray}

where

\begin{eqnarray}
\footnotesize
\begin{array}{cc}
\partial_k^{\beta(\alpha,[x_i]_k)}\delta_{jk}:= \der{\partial}{[x]_k}{\beta(\alpha,[x_i]_k)}\delta_{jk}, & 1\leq j,k\leq n,
\end{array}
\end{eqnarray}

with $ \delta_{jk} $ the Kronecker delta, $ \epsilon $ a positive constant $ \ll 1 $, and $ \beta (\alpha, [x_i]_k) $ a function defined as follows

\begin{eqnarray}\label{eq:c2.34}
\footnotesize
\begin{array}{cc}
\beta(\alpha,[x_i]_k):=\left\{
\begin{array}{cc}
\alpha, &\mbox{if \hspace{0.1cm} }  |[x_i]_k|\neq 0 \vspace{0.1cm}\\
1,& \mbox{if \hspace{0.1cm}  }  |[x_i]_k|=0
\end{array}\right..
\end{array}
\end{eqnarray}

It should be mentioned that the value $ \alpha = 1 $ in \eqref{eq:c2.34}, is taken to avoid the discontinuity that is generated when using the fractional derivative of constants in the value $ x = 0 $. Furthermore, since in the previous method  $\norm{\Phi^{(1)}(\alpha,\xi)}\neq 0$ if $\norm{f(\xi)}=0$, for the \textbf{Corollary \ref{cor:2-001}}, any sequence $ \set{x_i} _ {i = 0} ^ \infty $ generated by the iterative method \eqref {eq:c2.40} has an order of convergence (at least) linear. It should be noted that depending on the definition of fractional derivative used, the fractional iterative methods have the particularity that they may be used of local form \cite{gao2009local} or of global form \cite{torreshern2020}. These methods also have the peculiarity of being able to find complex roots of polynomials using real initial conditions \cite{fernando2017fractional}. Some differences between the classical Newton's method (CN), fractional Newton method (FN) \cite{torreshern2020} and fractional pseudo-Newton method (FPN) are listed in the Table \ref{tab:07}

\begin{table}[!ht]
\centering
\scriptsize
\begin{tabular}{c|c|c|c}
\toprule
&CN &
FN&
FPN \\ \midrule
\begin{tabular}{c}
Can it find \\
complex zeros\\
 of a polynomial \\
 using real initial\\
  conditions?
\end{tabular}& No & Yes & Yes \\ \midrule
\begin{tabular}{c}
Can it find\\
 multiple zeros \\
of a function \\
using a single \\
initial condition?
\end{tabular}& No & Yes & Yes \\ \midrule
\begin{tabular}{c}
Can it be \\
used if the\\
 function is not \\ differentiable?
\end{tabular}& No&Yes&Yes \\ \midrule
\begin{tabular}{c}
For a space of \\
dimension $N$.\\
How many (fractional)\\
partial 
 derivatives \\
 does it need?
\end{tabular}&
\begin{tabular}{c}
$N\times N$ 
\end{tabular}& 
\begin{tabular}{c}
$N\times N$
\end{tabular}
 & \begin{tabular}{c}
$N$ 
\end{tabular} \\ \midrule
\begin{tabular}{c}
Is it recommended \\
for solving systems\\
 where the (fractional)\\ partial derivatives\\  
 are analytically \\
 difficult to obtain?
\end{tabular}& No & No & Yes\\ \bottomrule
\end{tabular}
\caption{Some differences between the classical Newton's method and two fractional iterative methods.}\label{tab:07}
\end{table}

\section{Investment Under Uncertainty}

When a producer makes an investment to begin operating, it is important that he gets involved in how the use of technology can influence the expansion or reduction of his production, because expenses are incurred, whether he decides or not to do so, this is known as sunk cost. On the other hand, a producer can estimate how much he will be able to produce depending on prices, costs, and financial performance. But all of these values change over time; they are volatile. Therefore, producers face uncertain scenarios. The theory of investment under uncertainty scenarios \cite{dixit1994investment,price1999irreversible} proposes a model to determine when a producer must expand or close, depending on his income. The aforementioned model requires a series of equations shown below:

\begin{itemize}
\item[i)] Brownian movement of income:

\begin{eqnarray}
\footnotesize
\begin{array}{c}
dI=\left(\mu dt+ \sigma dz \right)I,
\end{array}
\end{eqnarray}

where

\begin{eqnarray*}
\footnotesize
\begin{array}{l}
I=\mbox{Income},\\
\mu=\mbox{Mean of continuous rate of income},\\
\sigma=
\mbox{Standard deviation of the continuous rate }\\
\hspace{0.5cm}\mbox{ of income movement}
,\\
dt=\mbox{Increase in time},\\
dz=\mbox{Variable that follows a Wiener process } \\
\hspace{0.6cm}\mbox{ with $N\to(0,1)$}.
\end{array}
\end{eqnarray*}

\item[ii)] Complementary function of the project value: 

\begin{eqnarray}
\footnotesize
\begin{array}{c}
V(I)=AI^{b} + BI^{a},
\end{array}
\end{eqnarray}

with

\begin{eqnarray*}
\footnotesize
\begin{array}{c}
b=-\dfrac{1}{\sigma^2}\left(    \left(\mu-\dfrac{1}{2}\sigma^2 \right)-\sqrt{\left(\mu-\dfrac{1}{2}\sigma^2 \right)^2+2\sigma^2 l  } \right), \vspace{0.1cm} \\
a=-\dfrac{1}{\sigma^2}\left(    \left(\mu-\dfrac{1}{2}\sigma^2 \right)+\sqrt{\left(\mu-\dfrac{1}{2}\sigma^2 \right)^2+2\sigma^2 l  } \right),
\end{array}
\end{eqnarray*}

where

\begin{eqnarray*}
\footnotesize
\begin{array}{l}
l=\mbox{Long-time real interest rate},\\
\end{array}
\end{eqnarray*}

it should be mentioned that $A$ and $B$ are parameters to be determined \cite{brambila2011bioeconomia}.

\item[iii)] The total function of the project value:

\begin{eqnarray}\label{eq:3-001}
\footnotesize
\begin{array}{c}
V(I)=AI^{b}+BI^{a}+\dfrac{I}{l-\mu}-\dfrac{c}{l},
\end{array}
\end{eqnarray}

where

\begin{eqnarray*}
\footnotesize
\begin{array}{c}
c=\mbox{Annual production cost}.
\end{array}
\end{eqnarray*}

\end{itemize}

Before continuing, it must be taken into account that there must be an income $H$, such that above this income it is decided to expand production and there must be an income $L$, such that below this income it is decided to reduce or close production. The value of the project when it is in $H$ or $L$ must fulfill the following conditions

\begin{eqnarray}\label{eq:3-002}
\footnotesize
\left\{
\begin{array}{c}
V_1(H)-V_0(H)-\kappa=0 \vspace{0.1cm} \\
V_1^{(1)}(H)-V_0^{(1)}(H)=0 \vspace{0.1cm}\\
V_1(L)-V_0(L)+\chi=0 \vspace{0.1cm}\\
V_1^{(1)}(L)-V_0^{(1)}(L)=0
\end{array}\right. ,
\end{eqnarray}

with

\begin{eqnarray}
\footnotesize
\begin{array}{c}
V_1^{(1)}(\ \cdot \ )-V_0^{(1)}( \ \cdot \ )=0,
\end{array}
\end{eqnarray}

the maximization condition between the final and initial cost function of the project, where

\begin{eqnarray*}
\footnotesize
\begin{array}{l}
\kappa=\mbox{Sunk cost},\\
\chi=\mbox{Cost of reducing or closing the project}.
\end{array}
\end{eqnarray*}

Finally, substituting the equation \eqref{eq:3-001} into the system \eqref{eq:3-002}, and considering the following notation

\begin{eqnarray*}
\footnotesize
\begin{array}{c}
(H,L,A,B)^T:=\left([x]_1,[x]_2,[x]_3,[x]_4\right)^T,
\end{array}
\end{eqnarray*}

the Dixit-Pindyck's system of equations \cite{dixit1994investment} is obtained

\begin{eqnarray}\label{eq:004}
\footnotesize
\left\{
\begin{array}{l}
a_5[x]_1+[x]_4[x]_1^{-a_1}-[x]_3[x]_1^{a_2}-a_6=0 \vspace{0.1cm} \\
-a_1[x]_4[x]_1^{-a_3}-a_2[x]_3[x]_1^{a_4}+a_5=0 \vspace{0.1cm} \\
a_5[x]_2+[x]_4[x]_2^{-a_1}-[x]_3[x]_2^{a_2}-a_7=0 \vspace{0.1cm}\\
-a_1[x]_4[x]_2^{-a_3}-a_2[x]_3[x]_2^{a_4}+a_5=0
\end{array}\right. ,
\end{eqnarray}

 which corresponds to the mathematical model of investment under uncertainty scenarios. This model was applied by Tauer for the case of milk producers \cite{tauer2004get}. The $ a_i $'s in the previous system are constants defined by the following expressions

 \begin{eqnarray*}
\left\{
\footnotesize
\begin{array}{c}
\begin{array}{cc}
a_1= \dfrac{\mu}{\sigma^2} -\dfrac{1}{2} + \rho ,&
a_2= -  \dfrac{\mu}{\sigma^2} +\dfrac{1}{2}+ \rho \vspace{0.1cm}\\
a_3= \dfrac{\mu}{\sigma^2} +\dfrac{1}{2} +  \rho,&
a_4= -  \dfrac{\mu}{\sigma^2} -\dfrac{1}{2}+ \rho \vspace{0.1cm}\\
\end{array}\\
\begin{array}{ccc}
a_5=\dfrac{1}{l-\mu},&a_6=\dfrac{c}{l}+\kappa,&
a_7=\dfrac{c}{l}-\chi
\end{array}
\end{array}
\right. ,
\end{eqnarray*}

with

\begin{eqnarray*}
\footnotesize
\begin{array}{c}
\rho=\sqrt{ \left(\dfrac{\mu}{\sigma^2}-\dfrac{1}{2}\right)^2+\dfrac{2l}{\sigma^2}   }.
\end{array}
\end{eqnarray*}

Through algebraic manipulations, the system \eqref{eq:004} may be rewritten as follows

\begin{eqnarray}\label{eq:006}
\scriptsize
\left\{
\begin{array}{l}
\left[x\right]_1
 =\dfrac{ a_6}{a_5}- \dfrac{a_1[x]_1^{a_2} \left( [x]_2^{a_3} - [x]_1^{a_3} \right) + a_2[x]_1 [x]_2^{a_3} \left([x]_1^{a_4} - [x]_2^{a_4}\right) }{a_1a_2 \left( [x]_1^{a_3 + a_4} - [x]_2^{a_3 + a_4} \right) } \vspace{0.1cm} \\
 \left[x\right]_2 =\dfrac{a_7}{a_5} -   \dfrac{a_1[x]_2^{a_2} \left( [x]_2^{a_3} - [x]_1^{a_3}\right) + a_2[x]_1^{a_3} [x]_2 \left([x]_1^{a_4} - [x]_2^{a_4}\right)}{a_1a_2\left([x]_1^{a_3 + a_4} - [x]_2^{a_3 + a_4}\right) }
\end{array},\right.
\end{eqnarray}

whose solution allows to know the values of the variables $[x]_3$ and $[x]_4$ through the following equations

\begin{eqnarray}\label{eq:0061}
\footnotesize
\left\{
\begin{array}{l}
\left[x\right]_3 = \dfrac{a_5 ([x]_1^{a_3} - [x]_2^{a_3})}{a_2 ([x]_1^{a_3 + a_4} - [x]_2^{a_3 + a_4})} \vspace{0.1cm}\\
\left[x\right]_4 = \dfrac{a_5 ([x]_1 [x]_2)^{a_3} ([x]_1^{a_4} - [x]_2^{a_4})}{a_1 ([x]_1^{a_3 + a_4} - [x]_2^{a_3 + a_4})}
\end{array}\right..
\end{eqnarray}

Using the system of equations \eqref{eq:006}, it is possible to define a  function $f:\Omega \subset \nset{R}^2\to \nset{R}^2$, that is,

\begin{eqnarray}\label{eq:007}
\scriptsize
\begin{array}{c}
f(x)=\begin{pmatrix}
\left[x\right]_1
 -\dfrac{ a_6}{a_5}+ \dfrac{a_1[x]_1^{a_2} \left( [x]_2^{a_3} - [x]_1^{a_3} \right) + a_2[x]_1 [x]_2^{a_3} \left([x]_1^{a_4} - [x]_2^{a_4}\right) }{a_1a_2 \left( [x]_1^{a_3 + a_4} - [x]_2^{a_3 + a_4} \right) } \vspace{0.1cm} \\
 \left[x\right]_2 -\dfrac{a_7}{a_5} +   \dfrac{a_1[x]_2^{a_2} \left( [x]_2^{a_3} - [x]_1^{a_3}\right) + a_2[x]_1^{a_3} [x]_2 \left([x]_1^{a_4} - [x]_2^{a_4}\right)}{a_1a_2\left([x]_1^{a_3 + a_4} - [x]_2^{a_3 + a_4}\right) }
\end{pmatrix}.
\end{array}
\end{eqnarray}

It should be mentioned that the system of equations \eqref{eq:006} depends in general on the values assigned to the constants $a_6$ and $a_7$, for this reason, for each new pair of assigned values it is necessary to calculate a new solution, that is,

\begin{eqnarray}
\footnotesize
\begin{array}{c}
(a_6,a_7)\overset{f}{\longrightarrow} x_n \in \nset{R}^2.
\end{array}
\end{eqnarray}

To solve the system \eqref{eq:006} using an iterative method, the task of finding a suitable initial condition must be carried out, which in many cases implies finding an initial condition near to the solution and this task usually takes some time. However, this problem is partially solved using the fractional pseudo-Newton method, because the initial condition does not need to be near to the solution \cite{torres2020approximation}, then the solution of the system may be determined more quickly.

\subsection{Some Results}

Assuming that system \eqref{eq:006} has the following constants

\begin{eqnarray*}
\footnotesize
\left\{
\begin{array}{l}
a_1=0.5355\\
a_2=1.5808\\
a_3=1.5355\\
a_4=0.5808\\
a_5=18.9753\\
\end{array}\right. ,
\end{eqnarray*}

and considering the following particular values

\begin{table}[!ht]
\centering
\footnotesize
$
\begin{array}{c}
\begin{array}{c|ccccc}
\toprule
&a_6&a_7&[x_0]_1&[x_0]_2 &||f(x_0)||_2\\ \midrule
1&451,474& 396,499 &15&20&6.00379e5\\
2&706,975 &  652,000 &17&18&9.61232e5\\
3&598,655 & 582,680 &9&16&8.35072e5 \\
4&506,975& 452,000 &5&19&6.78951e5\\ 
5&633,603&  578,628 &11&12&8.57733e5
\\
\bottomrule
\end{array} 
\end{array}
$
\caption{Different values of $a_6$ and $a_7$ along with some initial conditions $x_0$.}\label{tab:01}
\end{table}

We can choose the values from the Table \ref{tab:01} to use the iterative method given by  \eqref{eq:c2.40}. Consequently, we obtain the results of the Table \ref{tab:02}

\begin{table}[!ht]
\centering
\footnotesize
$
\begin{array}{c}
\begin{array}{c|ccc}
\toprule
&\alpha&[x_n]_1&[x_n]_2 \\ \midrule
1&0.26131&41,844.57090443&11,857.32126593\\
2&0.25628&60,324.4350877 &20,727.99532223\\
3&0.23116&43,561.70316013 &20,925.42239162\\
4&0.27136&45,951.77394332&13,741.03694719\\
5&0.24623&55,117.71562961 &18,133.15925118
 \\ \bottomrule
\end{array} 
\vspace{0.1cm}\\
\begin{array}{ccc}
\toprule
||x_n-x_{n-1}||_2&||f(x_n)||_2&n \\ \midrule
6.94046e-6	&6.96414e-5	&78\\
5.03627e-6&	8.61788e-5	&85\\
7.21678e-6&	8.66393e-5&	128\\
4.60353e-6	&8.71723e-5	&105\\
6.19923e-6	&9.26936e-5	&83
\\ \bottomrule
\end{array}
\end{array}
$
\caption{Results obtained using the iterative method \eqref{eq:c2.40} with $\epsilon=e-4$.}\label{tab:02}
\end{table}

\section{Conclusions}

The Dixit-Pindyck's system allows us to make the decision to continue or cancel the development of a project based on the observed income, which is essential when assuming the risk of making an investment under uncertain scenarios, but being a nonlinear system, iterative numerical methods are needed to find its solution and this generally implies that the task of finding an initial condition near to its solution must be done before trying to solve it. However, the fractional pseudo-Newton method has the characteristic that the initial condition does not need to be near to the solution to ensure convergence, added to the fact that it can find $N$ solutions of a system using a single initial condition, which makes it an ideal numerical method to solve the Dixit-Pindyck's system, which is the same as solving problems related to investment under uncertainty. Furthermore, this method, in contrast to Newton's method, is valid for non-differentiable functions and, consequently, can be used to model a greater number of phenomena related to physics, engineering and economics.

Partially funded by PAPPIT$\-$IT$\_$101421, UNAM.

\bibliography{Biblio}
\bibliographystyle{unsrt}

\end{document}